\documentclass[12pt,oneside]{amsart}

\usepackage[english]{babel}
\usepackage{hyphenat}
\usepackage[all]{xy}
\usepackage[utf8]{inputenc}
\usepackage[T1]{fontenc}
\usepackage{xspace}   

\usepackage{amsmath}
\usepackage{mathtools}
\allowdisplaybreaks
\usepackage{amssymb}
\usepackage{amscd}
\usepackage{amsthm} 
\usepackage{slashed}

\usepackage{graphicx} 

\graphicspath{{figures/}{/}}
\graphicspath{{figures/}{figuras/}} 

\usepackage{caption} 
\usepackage{subcaption}
\usepackage{wrapfig}
\usepackage{float}
\usepackage{epic,eepic}
\usepackage[table]{xcolor} 
\usepackage{tcolorbox}
\tcbuselibrary{listings,theorems}
\usepackage[absolute,overlay]{textpos}
\usepackage{MnSymbol}
\usepackage{tikz}
\usetikzlibrary{fit,patterns}
\tikzset{highlight/.style={rectangle,
                           fill=red!15,
                           blend mode = multiply,
                           rounded corners = 0.5 mm,
                           inner sep=1pt,
                           fit = #1}}

\usepackage{longtable}
\usepackage{tabularx}
\usepackage{tabularray}
\usepackage{multirow}
\usepackage{nicematrix} 

\usepackage[inline]{enumitem} 
\setlist{topsep=0pt} 
\usepackage{tasks}

\usepackage{url}
\usepackage{doi}
\usepackage{hyperref}
\hypersetup{
	urlcolor=blue,
}

\usepackage{fullpage}  
\usepackage[inline,final]{showlabels}

\usepackage[mathlines]{lineno} 

\theoremstyle{plain}
\newtheorem{thm}{Theorem}
\newtheorem{lem}[thm]{Lemma}
\newtheorem{cor}[thm]{Corollary}
\newtheorem{prop}[thm]{Proposition}

\newtheorem{rem}[thm]{Remark}
\newtheorem{claim}[thm]{Claim}

\newtheorem{obs}[thm]{Observation}

\newcommand{\dfn}[1]{{\em \bfseries #1}}

\newcommand{\ucr}[1]{\operatorname{cr}(#1)}
\newcommand{\rcr}[1]{\overline{\operatorname{cr}}(#1)}
\newcommand{\pcr}[1]{\widetilde{\operatorname{cr}}(#1)}
\newcommand{\symrcr}[1]{\operatorname{sym}-\overline{\operatorname{cr}}_3(#1)}
\newcommand{\sympcr}[1]{\operatorname{sym}-\widetilde{\operatorname{cr}}_3(#1)}

\newcommand{\pos}[1]{\operatorname{sp}(#1)} 
\newcommand{\pospi}[1]{\operatorname{sp}_{\pi}(#1)} 

\newcommand{\laposg}[1]{\overleftarrow{\operatorname{sp}_{>}}(#1)}
\newcommand{\laposl}[1]{\overleftarrow{\operatorname{sp}_{<}}(#1)}
\newcommand{\raposg}[1]{\overrightarrow{\operatorname{sp}_{>}}(#1)} 
\newcommand{\raposl}[1]{\overrightarrow{\operatorname{sp}_{<}}(#1)}

\newcommand{\set}[1]{\left\{#1\right\}}
\newcommand{\floor}[1]{\left\lfloor #1 \right\rfloor}
\newcommand{\ceil}[1]{\left\lceil #1 \right\rceil}

\newcommand{\real}{\mathbb{R}}

\title{The \(3\)-symmetric pseudolinear crossing number of \(K_{33}\)}

\author[V. H. G\'omez Mart\'inez]{V\'ictor H. G\'omez Mart\'inez}
\address{Facultad de Ciencias, Universidad Aut\'onoma de San Luis Potos\'i, San Luis Potos\'i, Mexico 78000.}
\thanks{Research supported by SECIHTI-M\'exico through a doctoral scholarship (Grant No. 370207)}
\email{a142080@alumnos.uaslp.mx}

\author[C. Hern\'andez-V\'elez]{C\'esar Hern\'andez-V\'elez}
\address{Facultad de Ciencias, Universidad Aut\'onoma de San Luis Potos\'i, San Luis Potos\'i, Mexico 78000.}
\email{cesar.velez@uaslp.mx}

\author[J. Lea\~nos]{Jes\'us Lea\~nos}
\address{Unidad Acad\'emica de Matem\'aticas, Universidad Aut\'onoma de Zacatecas, Zacatecas, Mexico 98000.}
\email{jleanos@uaz.edu.mx}

\begin{document}



\begin{abstract} We show that the $3$-symmetric rectilinear and the $3$-symmetric pseudolinear crossing numbers of $K_{33}$ are equal. Specifically, we prove that   \(\symrcr{K_{33}} = 14~634 = \sympcr{K_{33}}\). 
\end{abstract}
  
\maketitle

\section{Introduction}\label{s:intro}

The \dfn{crossing number} \(\ucr{G}\) of a graph \(G\) is the minimum number of pairwise crossing of edges in a drawing of \(G\) in the plane.   
The \dfn{rectilinear crossing number}  \(\rcr{G}\) of \(G\) is the minimum number of pairwise crossings of edges in a rectilinear drawing of \(G\) in the plane; that is, a drawing of \(G\) where the vertices of \(G\) are represented by points in general position, and an edge joining two vertices is represented by the straight segment joining the corresponding two points.

The problem of determining the exact value of~ \(\rcr{K_n}\) was proposed by 
Erd\H{o}s and Guy~\cite{erdos1973crossing} and remains open. Moreover, 
\(\rcr{K_n}\) is equal to the minimum number~ \(\square(n)\) of convex quadrilaterals determined by $n$ points in the plane in general position, which belongs to another well-known family of classical problems in combinatorial geometry, known as 
Erd\H{o}s–Szekeres type problems. \'Abrego et al.~\cite{abrego2012k} determine the exact value of \(\rcr{K_n}\) for \(n \leq 27\), and Cetina et al.~\cite{cetina2011point} for \(n = 30\). 
A table presenting lower bounds and best known minimizing examples for 
\(\rcr{K_n}\), with \(n\) up to 100, is shown in the Aichholzer's web page ``On the Rectilinear Crossing Number''~\cite{aichholzer-webpage}.

For the rest of the section, $P$ denotes an $n$-point set in the plane in general position. We recall that a {\em \(k\)-edge} of \(P\), with 
\(0 \leq k \leq (n-2)/2\), is a line through two points of \(P\) leaving exactly \(k\) points on one side. 
If \(n\) is even and \(k = (n-2)/2\), such \(k\)-edges are called \dfn{halving lines}. When \(n\) is odd the halving lines leave \((n-3)/2\) and \((n-1)/2\) points of \(P\) on each side.
A {\em \((\leq k)\)-edge} is any \(j\)-edge with \(0 \leq j \leq k\). As usual, we denote by \(E_{k}(P)\) the number of \(k\)-edges and by    \(E_{\leq k} (P)\) the number of (\(\leq k\))-edges of \(P\). We remark that \(E_{\leq \floor{n/2}}(P)=\binom{n}{2}\).

Let \(\rcr{P}\) be the number of crossings in the rectilinear drawing on \(K_n\) induced by \(P\).
Independently, \'Abrego and Fern\'andez-Merchant~\cite{abrego2005lower}, and Lov\'asz et al.~\cite{lovasz2004convex} found the following break-through connection between \(\rcr{P}\) and the number of \((\leq k)\)-edges:
\begin{align}
    \rcr{P} &= \sum_{k=0}^{\floor{\frac{n}{2}} - 2} (n - 2k - 3)E_{\leq k}(P) - \frac{3}{4} \binom{n}{3} + \left( 1 + (-1)^{n+1}\right) \frac{1}{8}\binom{n}{2},\label{eq:rcr-<=kedges}\\
    \intertext{or equivalently,}
    \rcr{P} &= 3\binom{n}{4} - \sum_{k=0}^{\floor{n/2}-1} k(n-k-2)E_{k}(P).\label{eq_rcr-kedges}
\end{align}
Let \(E_{\leq k} (n)\) denote the minimum of \(E_{\leq k} (P)\) taken over all \(n\)-point sets~\(P\). Determining \(E_{\leq k}(n)\) is another open problem in combinatorial geometry. The following lower bound for \(E_{\leq k}(n)\)  was independently established by \'Abrego and Fern\'andez-Merchant~\cite{abrego2005lower} and by Lov\'asz et al.~\cite{lovasz2004convex}: 
\begin{equation} \label{eq:k-egdes}
    E_{\leq k}(n) \geq 3 \binom{k + 2}{2}, \text{ for } 0 \leq k \leq n/2  - 2.
\end{equation}
Moreover, it is known that the lower bound in~\eqref{eq:k-egdes} is sharp 
for \(k < n/3\). Some improvements have been made for the lower bound of \(E_{\leq k}(n)\), 
in particular, Aichholzer et al.~\cite{aichholzer2007lower} proved that
\begin{equation}\label{eq:new_lower_k-edges}
    E_{\leq k}(n) \geq 3\binom{k+2}{2} + \sum_{j=\floor{n/3}}^k (3j - n +3).
\end{equation}

An equivalent form of inequality~\eqref{eq:new_lower_k-edges} was given by \'Abrego et al.~\cite{abrego2008extended}, in the more general context of simple generalized configurations of \(n\) points, showing that 
\begin{equation}\label{eq:lower-bound}
E_{\leq k} (n) \geq 3 \binom{k+2}{2} + 3 \binom{k +2 - \floor{n/3}}{2} - \max \left\{0, \left(k +1 - \floor{\frac{n}{3}}\right)\left(n - 3\floor{\frac{n}{3}}\right)\right\}. 
\end{equation}

Currently, the best lower bound is provided by Abrego et al.~\cite{abrego2012k}, whose results show that the bound in~\eqref{eq:lower-bound} is tight for \(k < \ceil{(4n - 11)/9}\). In~\cite{abrego2012k}, the authors exhibited sets of points that attain equality in~\eqref{eq:lower-bound} and proved that
\begin{equation}
E_{\leq k}(n) \geq \binom{n}{2} - \frac{1}{9}\sqrt{1 - \frac{2k + 2}{n}} \left(5n^2 + 19n - 31\right), \text{ for } k \geq \ceil{(4n-11)/9}.    
\end{equation}

In 2008, \'Abrego et al.~\cite{abrego2008maximum} established the following two bounds.
\begin{equation}\label{eq:halving}
E_{\floor{n/2}-1}(n) \leq 
\begin{cases}
    \floor{\frac{1}{2}\binom{n}{2} - \frac{1}{2} E_{\leq \floor{n/2}-3}(n)}, & \text{if } n \text{ is even},\\
    \floor{\frac{2}{3}\binom{n}{2} - \frac{2}{3} E_{\leq \floor{n/2}-3}(n) + \frac{1}{3}}, & \text{if } n \text{ is odd}.
\end{cases}
\end{equation}
\begin{equation}\label{eq:almost-halving}
E_{\leq \floor{n/2}-2}(n) \geq 
\begin{cases}
    \binom{n}{2} - \floor{\frac{1}{24}n(n+30)-3},& \text{if } n \text{ is even},\\
    \binom{n}{2} - \floor{\frac{1}{18}(n-3)(n+45)+ \frac{1}{9}},& \text{if } n \text{ is odd}.
\end{cases}
\end{equation}

By applying the lower bounds established in~\eqref{eq:lower-bound} and~\eqref{eq:almost-halving} within Equation~\eqref{eq:rcr-<=kedges} for a set of 33 points, we obtain \(\rcr{K_{33}} \geq 14~626\).
In contrast, Aichholzer provides on his web page~\cite{aichholzer-webpage} a configuration of 33 points for which the rectilinear crossing number of the induced complete graph is 14~634, which implies that 
\(14~626 \leq \rcr{K_{33}} \leq 14~634\).

All these parameters can be formulated in the more general context of the generalized configuration of points. 
We recall that a \dfn{pseudoline} is a simple curve in the proyective plane 
\(\mathbb{P}^2\) whose removal does not disconnect 
\(\mathbb{P}^2\).
An \dfn{arrangement of pseudolines} is a collection of pseudolines where every two of them intersect each other exactly once.
A \dfn{generalized configuration} is a set of points in the plane such that every pair of points has exactly one pseudoline passing through them, and the set of pseudolines forms an arrangement of pseudolines.

Let \(\mathcal{D}\) be a drawing of a graph \(G\) in the plane, and let \(C\) be a disk containing \(\mathcal{D}\). By discarding \(\real \setminus C\) and identifying antipodal points on the boundary of \(C\) we may consider 
\(\mathcal{D}\) as a drawing of \(G\) in \(\mathbb{P}^2\). If each edge of 
\(\mathcal{D}\) can be extended to a pseudoline so that the result is a pseudoline arrangement, we say that 
\(\mathcal{D}\) is a \dfn{pseudolinear drawing}. The \dfn{pseudolinear crossing number} \(\pcr{G}\) of \(G\) is the minimum number of pairwise crossing of edges in a pseudolinear drawing of \(G\). It is well known that a rectilinear drawing is pseudolinear. 
Because pseudolinear and rectilinear drawings are restricted classes of drawings, hence, for every graph \(G\) we have \(\ucr{G} \leq \pcr{G} \leq \rcr{G}\).

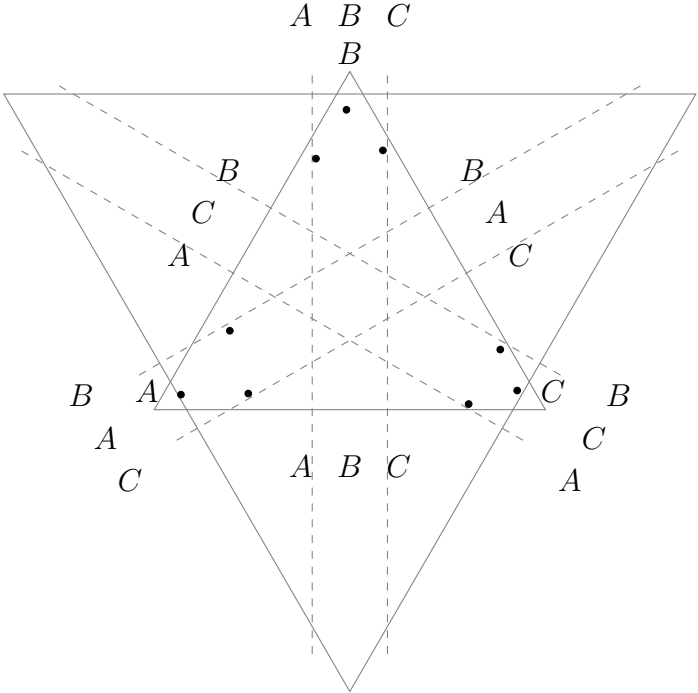
\begin{figure}[h]
\begin{center}

\begin{tikzpicture}[scale=0.5]
\draw[gray](0,6)--(-5.2,-3)--(5.2,-3)--(0,6);
\draw[gray](0,-10.5)--(-9.2,5.4)--(9.2,5.4)--(0,-10.5);
\draw[dashed,gray](-1,-9.5)--(-1,6);
\draw[dashed,rotate=120,gray](-1,-9.5)--(-1,6);
\draw[dashed,rotate=240,gray](-1,-9.5)--(-1,6);
\draw[dashed,gray](1,-9.5)--(1,6);
\draw[dashed,rotate=120,gray](1,-9.5)--(1,6);
\draw[dashed,rotate=240,gray](1,-9.5)--(1,6);
%
%
\draw(0,6.5) node{\(B\)};
\filldraw(-0.9,3.68) circle (2.5pt);
\filldraw(0.88,3.9) circle (2.5pt) ;
\filldraw(-0.09,4.98) circle (2.5pt);
\draw(5.4,-2.5) node{\(C\)};
\filldraw(4.45,-2.49) circle (2.5pt) ;
\filldraw(4,-1.4) circle (2.5pt) ;
\filldraw(3.16,-2.85) circle (2.5pt);
\draw(-5.4,-2.5) node{\(A\)};
\filldraw(-3.19,-0.9) circle (2.5pt) ;
\filldraw(-2.7,-2.57) circle (2.5pt) ;
\filldraw(-4.49,-2.6) circle (2.5pt);
\draw(-1.3,-4.5) node{\(A\)};
\draw(0,-4.5) node{\(B\)};
\draw(1.3,-4.5) node{\(C\)};
\draw[rotate=120](-1.3,-4.5) node{\(C\)};
\draw[rotate=120](0,-4.5) node{\(A\)};
\draw[rotate=120](1.3,-4.5) node{\(B\)};
\draw[rotate=240](-1.3,-4.5) node{\(B\)};
\draw[rotate=240](0,-4.5) node{\(C\)};
\draw[rotate=240](1.3,-4.5) node{\(A\)};
\draw(-1.3,7.5) node{\(A\)};
\draw(0,7.5) node{\(B\)};
\draw(1.3,7.5) node{\(C\)};
\draw[rotate=120](-1.3,7.5) node{\(C\)};
\draw[rotate=120](0,7.5) node{\(A\)};
\draw[rotate=120](1.3,7.5) node{\(B\)};
\draw[rotate=240](-1.3,7.5) node{\(B\)};
\draw[rotate=240](0,7.5) node{\(C\)};
\draw[rotate=240](1.3,7.5) node{\(A\)};
\end{tikzpicture}
\caption{A 9-points set which is \(3\)-decomposable.}
\label{fig:3symk12-3deck12}
\end{center}
\end{figure}

Although there is no conjectured value for \(\rcr{K_n}\), the structure of the drawings that are known to be optimal (crossing-minimal) is known.  
We say that a point set \(P\) is \dfn{3-decomposable} if it can be partitioned into three equal sized sets \(A, B\) and \(C\) such that there exists a triangle \(T\) enclosing \(P\) such that the orthogonal projection of \(P\) onto the three sides of \(T\) show \(A\) between \(B\) and \(C\) on one side, \(B\) between \(C\) and \(A\) on the second side, and \(C\) between \(A\) and \(B\) on the third side. See Figure~\ref{fig:3symk12-3deck12}. 

\begin{figure}[h]
\begin{center}
\begin{tikzpicture}
\begin{scope}
\draw[dashed,gray,rotate=15](0,0)--(0,1.8);
\draw[rotate=345](1.4,0.9) node{\(Q_{0}\)};
\filldraw[fill=white] (0.8,1) circle (1pt);
\filldraw[fill=white] (0.4,0.8) circle (1pt);
\filldraw[fill=white] (0.2,0.6) circle (1pt);
\filldraw[fill=white] (0.1,0.4) circle (1pt);
%
%
\draw[dashed,gray,rotate=135](0,0)--(0,1.8);
\draw[rotate=105](1.4,0.9) node{\(\rho^{2}(Q_{0})\)};
\filldraw[rotate=120] (0.8,1) circle (1pt);
\filldraw[rotate=120] (0.4,0.8) circle (1pt);
\filldraw[rotate=120] (0.2,0.6) circle (1pt);
\filldraw[rotate=120] (0.1,0.4) circle (1pt);
%
%
\draw[dashed,gray,rotate=255](0,0)--(0,1.8);
\draw[rotate=225](1.4,0.9) node{\(\rho(Q_{0})\)};
\filldraw[rotate=240] (0.8,1) circle (1pt);
\filldraw[rotate=240] (0.4,0.8) circle (1pt);
\filldraw[rotate=240] (0.2,0.6) circle (1pt);
\filldraw[rotate=240] (0.1,0.4) circle (1pt);
%
\end{scope}
\begin{scope}[shift={(7,0)},scale=1.5]
\draw[dashed,gray,rotate=180](-1,1)--(1,-1);
\filldraw[rotate=300](-0.7,-0.65) node {\(\rho(Q_{0})\)};
\filldraw[rotate=240](-0.7,-0.65) node {\(\rho^{2}(Q_{0})\)};
\filldraw[rotate=180](-0.7,-0.65) node {\(\rho^{3}(Q_{0})\)};
\filldraw[rotate=120](-0.7,-0.65) node {\(\rho^{4}(Q_{0})\)};
\filldraw[rotate=60](-0.7,-0.65) node {\(\rho^{5}(Q_{0})\)};
\draw(-0.7,-0.65) node {\(Q_{0}\)};

\filldraw (0.28,0.46) circle (0.75pt);
\filldraw (0.37,0.36) circle (0.75pt);
\filldraw (0.49,0.28) circle (0.75pt);
\filldraw (0.30,0.59) circle (0.75pt);

\draw[dashed,gray,rotate=60](-1,1)--(1,-1);
\draw[dashed,gray,rotate=120](-1,1)--(1,-1);
\filldraw[rotate=60] (0.28,0.46) circle (0.75pt);
\filldraw[rotate=60] (0.37,0.36) circle (0.75pt);
\filldraw[rotate=60] (0.49,0.28) circle (0.75pt);
\filldraw[rotate=60] (0.30,0.59) circle (0.75pt);
%
\filldraw[rotate=120] (0.28,0.46) circle (0.75pt);
\filldraw[rotate=120] (0.37,0.36) circle (0.75pt);
\filldraw[rotate=120] (0.49,0.28) circle (0.75pt);
\filldraw[rotate=120] (0.30,0.59) circle (0.75pt);
%
\filldraw[rotate=180,fill=white] (0.28,0.46) circle (0.75pt);
\filldraw[rotate=180,fill=white] (0.37,0.36) circle (0.75pt);
\filldraw[rotate=180,fill=white] (0.47,0.28) circle (0.75pt);
\filldraw[rotate=180,fill=white] (0.30,0.59) circle (0.75pt);
%
\filldraw[rotate=240] (0.28,0.46) circle (0.75pt);
\filldraw[rotate=240] (0.37,0.36) circle (0.75pt);
\filldraw[rotate=240] (0.49,0.28) circle (0.75pt);
\filldraw[rotate=240] (0.30,0.59) circle (0.75pt);
%
\filldraw[rotate=300] (0.28,0.46) circle (0.75pt);
\filldraw[rotate=300] (0.37,0.36) circle (0.75pt);
\filldraw[rotate=300] (0.49,0.28) circle (0.75pt);
\filldraw[rotate=300] (0.30,0.59) circle (0.75pt);
\end{scope}
\end{tikzpicture}
\caption{The underlying set of points of a \(3\)-symmetric drawing of \(K_{12}\) (left) and a \(6\)-symmetric drawing of \(K_{24}\) (right).}
\label{fig:6symkx-2symkx}
\end{center}
\end{figure}
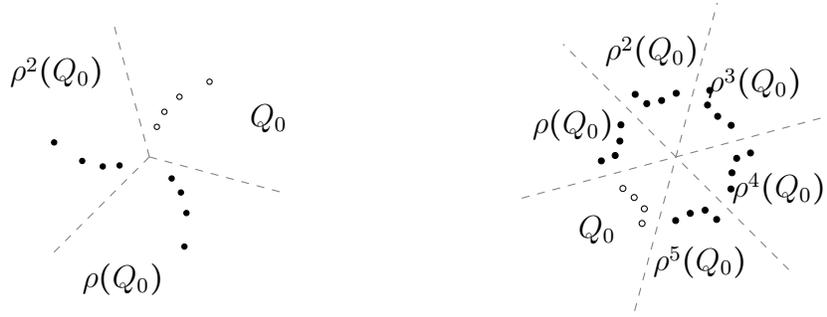

A point set \(P\) is \dfn{\(m\)-symmetric} if it can be partitioned into 
\(m\geq 2\) equal sized sets \(Q_0, Q_1, \ldots, Q_{m-1}\) invariant under a \(2\pi/m\)-rotation \(\rho\) (around a suitable point). Then, without loss of generality, we may assume that \(Q_i = \rho^{i}(Q_0)\), for
\(i = 0, \ldots, m-1\). We say that \(Q_0\) is the \dfn{seed set} of $P$. Clearly, if an \(n\)-point set is $3$-decomposable 
(\(m\)-symmetric, respectively), then \(n\) must be a multiple of $3$ ($m$, respectively). See Figure~\ref{fig:6symkx-2symkx}

A rectilinear drawing of \(K_n\) is 3-decomposable (\(3\)-symmetric, respectively) if its underlying point set is 3-decomposable (\(3\)-symmetric, respectively).

The exact value of \(\rcr{K_n}\) is known for very few values of \(n\). \'Abrego~\cite{abrego2008maximum} presented the exact value of \(\rcr{K_n}\), for \(n \leq 27\), and Cetina et al.~\cite{cetina2011point} for \(\rcr{K_{30}}\).
For every \(n\) multiple of 3, \(n \leq 30\), there is an optimal drawing of \(K_n\) that is 3-symmetric and 3-decomposable. Moreover, in general, for \(n\) multiple of 3, the best known (crossing-wise) drawings of \(K_n\) are 
3-decomposable. In fact, \'Abrego et al.~\cite{abrego20103sym} conjectured that, for each positive integer \(n\) multiple of 3, there is an optimal rectilinear drawing of \(K_n\) that is 3-symmetric and 3-decomposable.

Analogously, the \dfn{3-symmetric rectilinear crossing number} \(\symrcr{K_n}\) (\dfn{respectively, 3-symmetric pseudolinear crossing number} \(\sympcr{K_n}\)) of \(K_n\), with \(n\) multiple of 3, is the minimum number of pairwise crossing of edges in a 3-symmetric rectilinear (respectively, 3-symmetric pseudolinear) drawing of \(K_n\). It is also worth noting that, to date, there is no conjecture regarding the value of \(\symrcr{K_n}\). Recently, De \'Avila-Mart\'inez et al.~\cite{avila20243symmetric} proved that \(\sympcr{K_{36}} = \symrcr{K_{36}} = 21~174\).

In this paper, we determine the exact 3-symmetric pseudolinear crossing number of \(K_{33}\) and, as a consequence, its 3-symmetric rectilinear crossing number.

\begin{thm}\label{thm:main}
    \(\sympcr{K_{33}} = 14~634 = \symrcr{K_{33}}\).
\end{thm}

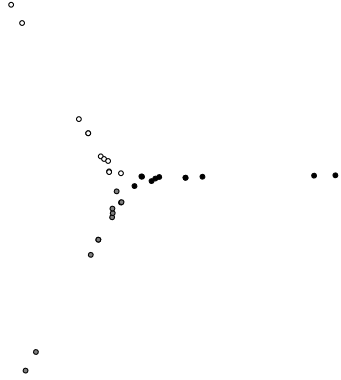
\begin{figure}[h]
\begin{center}

\begin{tikzpicture}[scale=0.03]
\filldraw[fill=white] (-50.0218885,79.3018474) circle (30pt);
\filldraw[fill=white] (-45.1723944,71.1948989) circle (30pt);
\filldraw[fill=white] (-20.0125330,28.5855310) circle (30pt);
\filldraw[fill=white] (-15.8721037,22.3132241) circle (30pt);
\filldraw[fill=white] (-15.8721040,22.31322453) circle (30pt);
\filldraw[fill=white] (-10.3183924,12.0586624) circle (30pt);
\filldraw[fill=white] (-8.8519236,10.9026774) circle (30pt);
\filldraw[fill=white] (-7.0502886,10.0103259) circle (30pt);
\filldraw[fill=white] (-6.6221918,5.3889958) circle (30pt);
\filldraw[fill=white] (-6.5940116,5.0836878) circle (30pt);
\filldraw[fill=white] (-1.3567216,4.5695226) circle (30pt);
%
\filldraw[rotate=120,fill=gray] (-50.0218885,79.3018474) circle (30pt);
\filldraw[rotate=120,fill=gray] (-45.1723944,71.1948989) circle (30pt);
\filldraw[rotate=120,fill=gray] (-20.0125330,28.5855310) circle (30pt);
\filldraw[rotate=120,fill=gray] (-15.8721037,22.3132241) circle (30pt);
\filldraw[rotate=120,fill=gray] (-15.8721040,22.31322453) circle (30pt);
\filldraw[rotate=120,fill=gray] (-10.3183924,12.0586624) circle (30pt);
\filldraw[rotate=120,fill=gray] (-8.8519236,10.9026774) circle (30pt);
\filldraw[rotate=120,fill=gray] (-7.0502886,10.0103259) circle (30pt);
\filldraw[rotate=120,fill=gray] (-6.6221918,5.3889958) circle (30pt);
\filldraw[rotate=120,fill=gray] (-6.5940116,5.0836878) circle (30pt);
\filldraw[rotate=120,fill=gray] (-1.3567216,4.5695226) circle (30pt);

\filldraw[rotate=240] (-50.0218885,79.3018474) circle (30pt);
\filldraw[rotate=240] (-45.1723944,71.1948989) circle (30pt);
\filldraw[rotate=240] (-20.0125330,28.5855310) circle (30pt);
\filldraw[rotate=240] (-15.8721037,22.3132241) circle (30pt);
\filldraw[rotate=240] (-15.8721040,22.31322453) circle (30pt);
\filldraw[rotate=240] (-10.3183924,12.0586624) circle (30pt);
\filldraw[rotate=240] (-8.8519236,10.9026774) circle (30pt);
\filldraw[rotate=240] (-7.0502886,10.0103259) circle (30pt);
\filldraw[rotate=240] (-6.6221918,5.3889958) circle (30pt);
\filldraw[rotate=240] (-6.5940116,5.0836878) circle (30pt);
\filldraw[rotate=240] (-1.3567216,4.5695226) circle (30pt);
\draw(-150,0)--(-150,0);
\draw(150,0)--(150,0);
\begin{scope}[shift={(-30,4)},every node/.style={scale=0.7}]
\draw(-223,70) node[right]{$Q=\{(-5002188850, 7930184740),$};
\draw(-200,58) node[right]{$(-4517239440,7119489890),$};
\draw(-200,46) node[right]{$(-2001253300,2858553100),$};
\draw(-200,34) node[right]{$(-1587210370,2231322410),$};
\draw(-200,22) node[right]{$(-1587210400,2231322453),$};
\draw(-200,10) node[right]{$(-1031839240,1205866240),$};
\draw(-200,-2) node[right]{$(-885192360,1090267740),$};
\draw(-200,-14) node[right]{$(-705028860,1001032590),$};
\draw(-200,-26) node[right]{$(-662219180,538899580),$};
\draw(-200,-38) node[right]{$(-659401160,508368780),$};
\draw(-200,-50) node[right]{$(-135672160,456952260)\}$};
\end{scope}
\end{tikzpicture}
\caption{The underlying set of points $P$ of an optimal 3-symmetric 3-decomposable rectilinear drawing of \(K_{33}\) with seed set \(Q\) and center of rotation at $(0,0)$. Namely, $\symrcr{P}=14~634$.}
\label{fig:K33}
\end{center}
\end{figure}

In Section~\ref{sec:allowable-seq} we recall the concept of allowable sequences and how the geometrical properties of a configuration of points are encoded by its corresponding allowable sequences. Hence, we analyze the allowable sequences of a pseudolinear 3–symmetric configuration of 33-points.
Section~\ref{sec:all-seq-min-K33} is devoted to characterize the allowable sequences that are optimal, that is, the induced generalized configuration is crossing-minimal.
An important property of the optimal allowable sequence is the 3-decomposability, which is studied in Section~\ref{sec:mo-bi-trans}.
Finally, in Section~\ref{sec:no-realizable} we prove that the best theoretically allowable sequence is not geometrically realizable by a 3-symmetric generalized configuration of 33-points, therefore, the 3-symmetric configuration of 33-points presented in this paper (see Figure~\ref{fig:K33}) is optimal.


\section{Allowable sequences}
\label{sec:allowable-seq}

An \(n\)-point set \(P\) in general position in the plane can be encoded by a \dfn{circular sequence} \(\boldsymbol{\Pi}_P\). 
As discussed in~\cite{lovasz2004convex}, we may assume without any loss of generality that the points in \(P\) are in \dfn{general position}; that is, no three points are collinear, and no two pairs of points lie on parallel lines.  We obtain a \dfn{halfperiod} \(\Pi\) of 
\(\boldsymbol{\Pi}_P\) as follows.
Let \(C\) be a circle containing \(P\) and let \(\ell\) be a directed tangent line to \(C\), and such that is not orthogonal to any line determined by two points in \(P\).
We project \(P\) orthogonally onto \(\ell\), and record the order of the point in \(P\)  on \(\ell\). 
This will be the first permutation \(\pi_0\) of \(\Pi\).
By rotating \(\ell\) clockwise, always keeping it tangent to \(C\), we follow the projection of \(P\) onto \(\ell\). 
Right after two points overlap in the projection, the order of \(P\) on \(\ell\) will change. 
This new order will be the permutation \(\pi_1\).
By continuing this process, until half a turn on \(C\), the order of the points in \(P\) on \(\ell\) will be the reverse of that in \(\pi_0\) and exactly \(\binom{n}{2}\) transpositions have taken place, one per each pair of points. 
The halfperiod \(\Pi\) (and, consequently, the circular sequence \(\boldsymbol{\Pi}\)) has the following properties:\\
\begin{enumerate*}[label=(SC\arabic*)]\\
\item\label{it:sc1} Consecutive permutations differ by a transposition of neighboring elements.\\
\item\label{it:sc2} The \(\binom{n}{2}\) permutations in 
\(\Pi\) are produced by the \(\binom{n}{2}\) possible transpositions of points in \(P\).
\end{enumerate*}

Let \( \boldsymbol{\Pi}_P = (\ldots, \pi_{-1}, \pi_0, \pi_1, \ldots) \) be the circular sequence associated to \(P\).
Let \(\tau_j\) be the transposition of adjacent elements, say \(x\) and \(y\), that transforms \(\pi_{j-1}\) into \(\pi_j\).
We write \(\tau_j = x|y\) to indicate the transposition of \(x\) and \(y\) where \(x\) originally precedes \(y\) (thus, \(x|y\) and \(y|x\) are different transpositions).
The properties~\ref{it:sc1} and~\ref{it:sc2} imply that when transposition \(x|y\) occurs then after \(\binom{n}{2}\) permutations the transposition \(y|x\) must occur. Then 
\(\pi_j, \pi_{j + \binom{n}{2}} \in \boldsymbol{\Pi}\) are reverse permutations of each other; in particular, if two elements transpose, they do not transpose again until every other pair has been transposed.

Goodman and Pollack~\cite{goodman1980combinatorial} defined \dfn{allowable sequence} as an infinite sequence of permutations of \(n\) elements that satisfies both properties~\ref{it:sc1} and~\ref{it:sc2}; that is, an allowable sequence is the combinatorial abstraction of a circular sequence. Moreover, an allowable sequence is determined by its \dfn{transposition sequence} \(T(\boldsymbol{\Pi}) = (\ldots, \tau_{-1}, \tau_0, \tau_1, \ldots)\), where \(\tau_j\) is the transposition of neighboring elements that transforms \(\pi_{j-1}\) into 
\(\pi_j\). We remark that an allowable sequence is completely determined by \(\binom{n}{2}\) consecutive permutations, which is called a \dfn{halfperiod}. In~\cite{goodman1980combinatorial}, it was also proved that there exist allowable sequences that are not circular sequences, and that there is a one-to-one correspondence between allowable sequences and generalized configurations of points.

Let \(\boldsymbol{\Pi}\) be a circular sequence on $n\geq 3$ elements, and let $i\in \{1,2,\ldots , n-1\}$. A transposition that occurs between elements in positions \(i\) and \(i + 1\) is an \dfn{\(i\)-transposition}. If \(i \leq n/2\), then an \dfn{\(i\)-critical transposition} is either an \(i\)-transposition or an \((n-i)\)-transposition, and a \dfn{\((\leq k)\)-critical transposition} is any \(i\)-critical transposition with \(i \leq k\).
If \(\Pi\) is a halfperiod of \(\boldsymbol{\Pi}\), then 
\(N_k(\Pi)\) denotes the number of \(k\)-critical transposition and \(N_{\leq k}(\Pi)\) denotes the number of \((\leq k)\)-critical transpositions in \(\Pi\), respectively.

\begin{obs}
Let \(P\) be an \(n\)-point set in the plane in general position and let \(1\leq k < n/2\). If $\Pi_P$ is a halfperiod of the circular sequence associated to $P$, then
\[E_{\leq (k-1)}(P)  = N_{\leq k} (\Pi_P).\]
\end{obs}


These parameters can be extended to generalized configurations of points, or equivalently, to allowable sequences \(\boldsymbol{\Pi}\).
In this context, an \(i\)-critical transposition corresponds to a \((k-1)\)-edge.
Therefore, \(E_{\leq k}(\boldsymbol{\Pi})\) corresponds to the number of \((k+1)\)-critical transpositions in any halfperiod of \(\boldsymbol{\Pi}\). Also, the concepts of 3-decomposability and 3-symmetry have been extended to the setting of allowable sequences~\cite{abrego2012k, avila20243symmetric}. 

An allowable sequence \(\boldsymbol{\Pi}\) on \(n\)-elements, \(n\) multiple of 3, is \dfn{3-decomposable} if there exists a halfperiod \(\Pi\) of \(\boldsymbol{\Pi}\) such that the elemets of \(\Pi\) can be labeled \(A = \{a_1, a_2, \ldots, a_{n/3}\}\), \(B = \{b_1,b_2, \ldots, b_{n/3}\}\), and \(C = \{c_1, c_2, \ldots, c_{n/3}\}\) so that \[\pi_0 = (a_1, a_2, \ldots, a_{n/3}, b_1, b_2, \ldots, b_{n/3}, c_1, c_2, \ldots, c_{n/3}),\] and there are \(0 < i < j \leq \binom{n}{2}\) such that 
\(\pi_{i+1}\) shows all the \(b\)-elements followed by all the \(a\)-elements followed by all the \(c\)-elements, and \(\pi_{j+1}\)  shows all \(b\)-elements followed by all the \(c\)-elements followed by all the \(a\)-elements.

Similarly, an allowable sequence \(\boldsymbol{\Pi}\) on an \(n\)-set $S$, \(n\) multiple of 3, is \dfn{3-symmetric} if there exists a halfperiod \(\Pi\) of \(\boldsymbol{\Pi}\), a subset \(U\subset S\) with \(n/3\) elements, and a permutation 
\(\sigma \colon\, S \to S\) such that:
\begin{enumerate}[topsep=0pt]
    \item \(\sigma^3\) is the identity map on \(S\),
    \item \(S\) is the disjoint union of \(U, \sigma(U)\) and \(\sigma^2(U)\),
    \item If \(x\) and \(y\) are distinct elements of \(S\) and \(x|y\) is a \(k\)-critical transposition, then 
    \(\sigma(x)|\sigma(y)\) is also a \(k\)-critical transposition.
    \item \(\sigma\) preserves the cyclic order in which any element of \(x \in S\) transposes with the elements in \(S \setminus \{x\}\); that is, if \(x \in S\) transposes with \(x_1, x_2, \ldots, x_{n-1} \in S \setminus \{x\}\) in that given order, then \(\sigma(x)\) transposes with \(\sigma(x_{t+1}), \sigma(x_{t+2}),\ldots, \sigma(x_{n-1}), \sigma(x_1), \ldots, \sigma(x_t) \in S \setminus \{\sigma(x)\}\) in that given order, for some integer \(0 \leq t \leq n-1\). 
\end{enumerate}

\begin{obs}
If \(\boldsymbol{\Pi}\) is a 3-symmetric allowable sequence on an \(n\)-set, then the number \(N_{\leq k}(\boldsymbol{\Pi)}\) of \((\leq k)\)-critical transpositions is a multiple of 3.
\end{obs}

Let \(\Pi=(\pi_0, \pi_{1}, \ldots , \pi_{\binom{n}{2}})\) be a fixed \(n\)-halfperiod of \({\bf \Pi}\), and let \(x\) and \(y\) be distinct elements of~\(P\). 
For \(\pi_j\in \Pi\), we will use \(\pi_j(x)\) to denote the position of \(x\) in \(\pi_j\); then, \(\pi_j(x)\in \{1,2,\ldots ,n\}\).
It follows that if \(\pi_0(x)<\pi_0(y)\), then there exists a unique ordered pair 
\[(j,i)\in \left\{1,2,\ldots , \binom{n}{2}\right\}\times \{1,2,\ldots ,n-1\},\] 
such that
\begin{enumerate}[topsep=0pt]
\item \(\pi_{j-1}(x)=i\) and \(\pi_{j-1}(y)=i+1\), 
\item \(\pi_{j}(x)=i+1\) and \(\pi_{j}(y)=i\), and
\item \(\pi_{j-1}(z)=\pi_{j}(z)\) for any \(z\in P \setminus\{x,y\}\). 
\end{enumerate}

We recall that \(\tau_{j}\) is the transposition that transforms \(\pi_{j-1}\) into \(\pi_{j}\). To emphasize the dependence on \(x\) and \(y\), we often write 
\(x|y\) instead of \(\tau_j\) and say that \(x|y\) occurred in the \dfn{\(i\)-gate}, which is abbreviated by \(g(x|y)=i\).
Then, the set of \(k\)-critical transpositions of \(\Pi\) is equal to \(\{x|y \ \colon\  g(x|y) \in \{k,n-k\}\}\).


\section{On allowable sequences that minimize the crossing number of \texorpdfstring{\(K_{33}\)}{K33}}
\label{sec:all-seq-min-K33}

For the remainder of this paper, any halfperiod under considerations belongs to an allowable sequence on a 
$33$-set $P$. We recall that
\[14~626 \leq \pcr{K_{33}} \leq \rcr{K_{33}} \leq 14~634.\]
Since there exists a 3-symmetric rectilinear drawing of \(K_{33}\) with 14~634 crossing (namely, that in Figure~\ref{fig:K33}), then \[14~626 \leq \pcr{K_{33}} \leq \sympcr{K_{33}} \leq \symrcr{K_{33}} \leq 14~634.\]

The \dfn{\((\leq k)\)-transpositions vector} of  \(\Pi\) is the vector \(\boldsymbol{V}_{\leq k}(\Pi) \coloneqq (N_{\leq 1}(\Pi), N_{\leq 2}(\Pi),\ldots, N_{\leq 16}(\Pi))\).
The \dfn{\(k\)-transpositions vector} of \(\Pi\) is the vector \(\boldsymbol{V}_k(\Pi) \coloneqq (N_{1}(\Pi), N_{2}(\Pi),\ldots, N_{16}(\Pi))\).
Let \(\Pi\) and \(\Pi'\) be two halfperiods. We write \(\boldsymbol{V}_{\leq k}(\Pi) \leq \boldsymbol{V}_{\leq k}(\Pi')\) if \(N_{\leq j}(\Pi) \leq N_{\leq j}(\Pi')\) for every \(j \in \set{1,2,\ldots,16}\).

From bounds~\eqref{eq:lower-bound} and~\eqref{eq:almost-halving}, we know that any halfperiod \(\Pi\) satisfies the following:
\[\boldsymbol{V}_{\leq k}(\Pi) \geq (3,9,18,30,45,63,84,108,135,165,198,237,282,333,398,453).\]

The 15th entry of this vector is given by~\eqref{eq:almost-halving}, while the remaining entries are provided by~\eqref{eq:lower-bound}.

If \(\Pi\) is a 3-symmetric halfperiod, then \(N_{\leq k}(\Pi)\) must be a multiple of 3 for every \(k  \in \set{1,2,\ldots, 16}\). In particular, \(N_{\leq 15}(\Pi)\) must be a multiple of 3, so \(N_{\leq 15}(\Pi) \geq 399\). This last fact and the identity 
 $N_{\leq 16}(\Pi)=528$ imply the following.

\begin{rem}\label{rem:vector-3symm}
If $\Pi$ is 3-symmetric halfperiod of an allowable sequence on a $33$-set, then
\begin{align*}
\boldsymbol{V}_{\leq k}(\Pi) &\geq (3,9,18,30,45,63,84,108,135,165,198,237,282,333,\boldsymbol{399},528).\\
\end{align*}
\end{rem}

By substituting the values of $N_{\leq k}(\Pi)$ provided in Remark~\ref{rem:vector-3symm} into the following equation (which is equivalent to Equation~\ref{eq:rcr-<=kedges}), we obtain:
\begin{equation}
    \pcr{\Pi} = \sum_{k=1}^{\floor{\frac{n}{2}}-1} (n - 2k - 1)N_{\leq k}(P) - \frac{3}{4} \binom{n}{3} + \left( 1 + (-1)^{n+1}\right) \frac{1}{8}\binom{n}{2}=14~628.\label{eq:pcr-<=ktrans}
\end{equation}
Then,
\[\boldsymbol{14~628} \leq \sympcr{K_{33}} \leq \symrcr{K_{33}} \leq 14~634.\]

\begin{prop}\label{prop:<=k-egdes_tight}
If $\Pi$ is a 3-symmetric halfperiod  of an allowable sequence on a $33$-set such that $\pcr{\Pi}<14~634$, then 
\begin{align*}
\boldsymbol{V}_{\leq k}(\Pi) &= (\underset{k=1}{3},\underset{k=2}{9},\underset{k=3}{18},\underset{k=4}{30},\underset{k=5}{45},\underset{k=6}{63},\underset{k=7}{84},\underset{k=8}{108},\underset{k=9}{135},\underset{k=10}{165},\underset{k=11}{198},\underset{k=12}{237},\underset{k=13}{282},\underset{k=14}{333}, \underset{k=15}{399},\underset{k=16}{528})\label{vec_trans1}\tag{V1}.\\ 
\boldsymbol{V}_{k}(\Pi) &= (\underset{k=1}{3},\underset{k=2}{6},\underset{k=3}{9},\underset{k=4}{12},\underset{k=5}{15},\underset{k=6}{18},\underset{k=7}{21},\underset{k=8}{24},\underset{k=9}{27},\underset{k=10}{30},\underset{k=11}{33},\underset{k=12}{39},\underset{k=13}{45},\underset{k=14}{51},\underset{k=15}{66},\underset{k=16}{129})\label{vec_trans2}\tag{V2}.
\end{align*}
\end{prop}
\begin{proof} Clearly, $N_{\leq 16}(\Pi)=528$. Since $\Pi$ is 
3-symmetric, each $N_{\leq k}(\Pi)$ must be a multiple of 3 for all
$k$. A simple calculation using 
Equation~\ref{eq:pcr-<=ktrans} indicates that if, for some 
\(k \in \{1,2, \ldots ,14\}\), the value of
$N_{\leq k}(\Pi)$ exceeds the claimed amount, then 
\(\sympcr{\Pi} \geq 14~640\). 
Similarly, from Equation~\ref{eq:pcr-<=ktrans}, it follows that if 
\(N_{\leq 15}(\Pi) \geq 402\), then \(\sympcr{\Pi} \geq 14~634\).

The second vector follows from the first one and consider that 
\(N_{k}(\Pi)=N_{\leq k}(\Pi)- N_{\leq k-1}(\Pi)\) for each 
$k\in \{1,2,\ldots, 16\}$.
\end{proof}

Cetina et al.~\cite{cetina2011point} proved the following relationship between 3-decomposability and the number of \((\leq k)\)-edges (equivalently \((\leq (k+1))\)-critical transpositions).

\begin{thm}[Main Theorem, Cetina et al.~\cite{cetina2011point}]\label{thm:<=k-edges_3decomp}
Let \(P\) be an \(n\)-point set with \(n\) multiple of 3. If \(P\) has exactly \(3\binom{k+2}{2}\) \((\leq k)\)-edges, for all \(0 \leq k < n/3\), then \(P\) is 3-decomposable.    
\end{thm}

Proposition~\ref{prop:<=k-egdes_tight} and Theorem~\ref{thm:<=k-edges_3decomp} imply the following.

\begin{obs}\label{obs:k33_opt-3decomp}
If $\Pi$ is a 3-symmetric halfperiod  of an allowable sequence on a $33$-set such that $\pcr{\Pi}<14~634$, then $\Pi$ is 3-decomposable.  
\end{obs}

\begin{rem}\label{rem:Pi}
The rectilinear drawing induced by the 33-point set in Figure~\ref{fig:K33} shows that $\sympcr{K_{33}}\leq 14~634$. 
Given this fact and Observation~\ref{obs:k33_opt-3decomp}, 
we will assume for the remainder of the paper that 
\(\Pi = (\pi_{0},\pi_{1},\ldots , \pi_{\binom{33}{2}})\) is a 
3-symmetric 3-decomposable halfperiod of an allowable sequence on a $33$-set $P$ such that $\sympcr{K_{33}}<14~634$.    
\end{rem}


\section{Monochromatic and bichromatic transpositions}
\label{sec:mo-bi-trans}

Throughout this section, \(\Pi\) and $P$ are defined as in Remark~\ref{rem:Pi}. Assume that \(A=\{a_{11},\ldots ,a_{1}\},\-B=\{b_{1},\ldots ,b_{11}\}\) and 
\(C=\{c_{1},\ldots ,c_{11}\}\) form a $3$-decomposition of $\Pi$. 

Let $\tau$ be a transposition of $\Pi$. We say that $\tau$ is \dfn{monochromatic} if it involves two elements from 
$X\in \{A,B, C\}$, otherwise $\tau$ is \dfn{bichromatic}. We denote the number of monochromatic (bichromatic, respectively) \((\leq k)\)-critical transpositions of \(\Pi \) by \(N_{\leq k}^{\rm mo} (\Pi )\), (\(N_{\leq k}^{\rm bi} (\Pi )\), respectively). 

\begin{rem}\label{r:N_k}
From the involved definitions, it is easy to see that
\begin{align*}
    N_{\leq k} (\Pi ) &= N_{\leq k}^{\rm mo}(\Pi) + N_{\leq k}^{\rm bi}(\Pi ),\\
    N_{k}^{\rm bi}(\Pi) &= N_{\leq k}^{\rm bi}(\Pi) - N_{\leq (k-1)}^{\rm bi}(\Pi),\\
    N_{k}^{\rm mo}(\Pi) &= N_{\leq k}^{\rm mo}(\Pi) - N_{\leq (k-1)}^{\rm mo}(\Pi).
\end{align*}
\end{rem}

\'Abrego et al.~\cite{abrego20103sym} proved the following.

\begin{prop}[Claim~1, \'Abrego et al.~\cite{abrego20103sym}]\label{prop:bichromatic}
Let \(\Pi^{'}\) be a 3-decomposable halfperiod on \(n\)~points, \(n\)  multiple of 3, and let \(k < n/2\). Then
\begin{equation}\label{eq:bichromatics}
N_{\leq k}^{\rm bi}(\Pi^{'}) =
\begin{cases}
3\binom{k+1}{2} &\text{if } k \leq n/3,\\
3\binom{n/3 + 1}{2} + \left(k - n/3\right) n&\text{if } n/3 < k < n/2.    
\end{cases}
\end{equation}
\end{prop}
From previous proposition, we have the following. 
\begin{cor}\label{cor:bichromatics-k33} Let $\Pi$ be as in Remark~\ref{rem:Pi}. Then, 
\begin{equation}\label{eq:bichromatic-k33}
N_{k}^{\rm bi}(\Pi) =  
\begin{cases}
    3k & \text{for } k \in \{1, \ldots , 11\}\\
    33 & \text{for } k \in \{12, \ldots , 16\}.
\end{cases}
\end{equation}
\hfill \qedsymbol{}
\end{cor}
The following corollary is a direct consequence of Proposition~\ref{prop:<=k-egdes_tight}, Remark~\ref{r:N_k} and Corollary~\ref{cor:bichromatics-k33}.

\begin{cor}\label{cor:monochromatics2-4-6} Let $\Pi$ be as in Remark~\ref{rem:Pi}. Then,  
\(N_k^{\rm mo}(\Pi) = 0\), for \(k\in\{1, \ldots, 11\}\), 
\(N_{12}^{\rm mo}(\Pi) = 6\), \(N_{13}^{\rm mo}(\Pi) = 12\), \(N_{14}^{\rm mo}(\Pi) = 18\), \(N_{15}^{\rm mo}(\Pi) = 33\), and \(N_{16}^{\rm mo}(\Pi)=96\).
\end{cor}

We will denote by $N_{k}^{aa}(\Pi)$ (respectively, $N_{\leq k}^{aa}(\Pi)$) the number of monochromatic $k$-critical 
($(\leq k)$-critical, respectively) transpositions that occur between elements of 
$A$. Analogously, we define $N_{k}^{bb}(\Pi), N_{k}^{cc}(\Pi), N_{\leq k}^{bb}(\Pi)$ and $N_{\leq k}^{cc}(\Pi)$. Then,
\(N_{k}^{\rm mo}(\Pi) = N_{k}^{aa}(\Pi) + N_{k}^{bb} + N_{k}^{cc}(\Pi)\) and \(N_{\leq k}^{\rm mo}(\Pi) = N_{\leq k}^{aa}(\Pi) + N_{\leq k}^{bb} + N_{\leq k}^{cc}(\Pi)\).

The following statements are straightforward consequences of Proposition~\ref{prop:<=k-egdes_tight}, Corollaries~\ref{cor:bichromatics-k33}, \ref{cor:monochromatics2-4-6},
and the $3$-symmetry of $\Pi$.  

\begin{cor}\label{cor:xx}
If $\Pi$ is as in Remark~\ref{rem:Pi}, then 
  \(N_{12}^{xx}(\Pi)~=~2\), \(N_{13}^{xx}(\Pi)~=~4\), \(N_{14}^{xx}(\Pi)~=~6\), \(N_{15}^{xx}(\Pi)~=~11\), \(N_{16}^{xx}(\Pi)~=~32\), for \(x \in \set{a,b,c}\).
\end{cor}


\section{There is no halfperiod satisfying conditions of Remark~\ref{rem:Pi}: Proof of Theorem~\ref{thm:main}}
\label{sec:no-realizable}

 In this section, \(\Pi\) and $P$ are defined as in Remark~\ref{rem:Pi}. Let 
 $\{A,B,C\}$ be a $3$–decomposition of $\Pi$. We assume that in the initial permutation $\pi_0$ of $\Pi$, the set \(B\) is between \(A\) and \(C\) (i.e., \(\pi_0 = (A, B, C)\)). Furthermore, we may assume that any transposition in $\Pi$ between an element of $A$ and an element of $B$ occurs before any transposition between an element of $A$ and an element of $C$, and the last bichromatic transpositions are those between an element in $B$ and an element in $C$. We denote this order by \(AB \prec AC \prec BC\).

Let $x|y$ be a transposition of $\Pi$. If $g(x|y)\leq 11$ (respectively, $g(x|y)\geq 22$), we say that $x|y$ occurs in the \dfn{initial third} (respectively, \dfn{final third}). Otherwise, we say that $x|y$ occurs in the \dfn{middle third}. According to Proposition~\ref{cor:monochromatics2-4-6}, we know that 
\(N_{k}^{\rm mo}(\Pi) = 0\) for \(k \in \{1, \ldots ,11\}\). Thus, every monochromatic transposition occurs in the middle third; therefore, in the initial and final third only bichromatic transpositions can occur.
These facts imply that exactly one element, which we label $x_{11}$, of 
$X\in \{A,B,C\}$ can occupy the $1$-st or $33$-th positions in any permutation of $\Pi$. Similarly, for each $i\in \{2,3,\ldots ,11\}$, there is a unique element of $X$ that reaches the positions $i$ or $34-i$ of some permutation of $\Pi$ but in none of them occupies the positions $i-1$ nor $35-i$. 
Let $x_{12-i}$ denote such an element of $X$. Then, 
\(A=\{a_{1},\ldots ,a_{11}\},~B=\{b_{1},\ldots ,b_{11}\}\), 
~\(C=\{c_{1},\ldots ,c_{11}\}\) and $\pi_0$ looks like:  
\[\pi_{0}=(a_{11},a_{10},\ldots ,a_{1},b_{l_1},b_{l_2},\ldots ,b_{l_{11}}\-,c_{1},c_{2},\ldots c_{11}),\] where \(j\in\{1,2, \ldots, 11\}\) and $b_{l_j}$ denotes the element of $B$ that is in the position $11+j$ of $\pi_0$.

For \(i\in\{1,2, \ldots, 11\}\), we say that an element \(x\) \dfn{exits} through the \(i\)th gate if it moves from position \(i\) to \(i + 1\); and for \(i \in \{23,24, \ldots ,33\}\), \(x\) \dfn{exits} through the \(i\)th gate if it moves from position  \(i\) to \(i - 1\).   We say that \(x\) \dfn{enters the middle third} if it exits through the 11th gate or through the 23rd gate.
Similarly, we say that \(x\) \dfn{leaves the middle third} if it moves from position 12 to 11, or from position 22 to 23.

\begin{prop}\label{prop:inc-dec}
For $X\in \{A,B,C\}$, the following hold:
\begin{itemize}
\item[(1)] The elements of $X$ enter (respectively, leave) the middle third in increasing (respectively, decreasing) order.
\item[(2)] $g(a_{11-r}|b_{r+1})=11$ for $r\in \{0,1,2, \ldots, 10\}$.
\item[(3)] $g(a_{11-r}|c_{r+1})=22$ for $r\in \{0,1,2, \ldots, 10\}$. 
\item[(4)] $g(b_{11-r}|c_{r+1})=11$ for $r\in \{0,1,2, \ldots, 10\}$. 
\end{itemize}
\end{prop}
\begin{proof}
Assertion (1) follows from the absence of monochromatic transpositions in the initial and final third, as well as from the labeling of the elements in $A\cup B \cup C$. 

For $r\in \{0,1,2, \ldots, 10\}$, let $\tau_{j_{r}}$ be the transposition in which $a_{11-r}$ enters the middle third. Clearly, $g(\tau_{j_r})=11$.  
Since $\tau_{j_r}$ does not occur in the middle third, it must be 
bichromatic. From \(AB \prec AC \prec BC\) it follows that $\tau_{j_r}=a_{r+1}|b_k$ for some $k\in \{1,\ldots, 11\}$. From (1) 
and the absence of monochromatic transpositions in the initial third
it follows that $k=11-r$. Assertions (3) and (4) follow in a similar manner.
\end{proof}

\begin{rem}\label{rem:same-a-b-c}
The 3-decomposability and the 3-symmetry of $\Pi$, along with the labeling of the elements in $A\cup B \cup C$ imply that everything that we show for $a_i|a_j$ also holds for $b_i|b_j$ and  $c_i|c_j$.
\end{rem}

By convention, if $j$ and $l$ are nonnegative integers such that 
$l<j$, then the set $\{v_j,\ldots, v_l\}$ will be the empty set. 

\begin{prop}\label{eta_12+r}
Let \(\Pi\) be  as in Remark~\ref{rem:Pi}. The following hold:
\begin{itemize}
\item[(1)] $\Pi$ has (exactly) 2 transpositions that move $a_{11}$ from left to right and contribute to \(N_{k}^{aa}(\Pi)\) for each $k\in \{12,13,\ldots ,16\}$.  
\item[(2)] $\Pi$ has (exactly) 2 transpositions that move $a_{10}$ from left to right and contribute to \(N_{k}^{aa}(\Pi)\) for each $k\in \{13,14\}$.  Additionally, $g(a_{10}|b_{1})=12$ and $g(a_{10}|c_{1})=21$.
\item[(3)] $\Pi$ has (exactly) 2 transpositions that move $a_{9}$ from left to right and contribute to \(N_{14}^{aa}(\Pi)\). Additionally, $g(a_{9}|b_{2})=12$,  $g(a_{9}|b_{1})=13$, $g(a_{9}|c_{1})=20$, $g(a_{9}|c_{2})=21$.
\item[(4)] For \(r\in\{0,1,2\}\), any transposition of \(\Pi\) contributing to 
\(N_{12+r}^{aa}(\Pi)\) involves exactly one element of \(\{a_{11-r},\ldots ,a_{11}\}\).  
\end{itemize}
\end{prop}
\begin{proof} Let $r\in \{0,1,2\}$. From the first two assertions of Proposition~\ref{prop:inc-dec}, we know that $a_{11-r}$ enters (respectively, leaves) the middle third with $\tau_{r_e}:=a_{11-r}|b_{r+1}$ (respectively, 
$\tau_{r_l}:=a_{11-r}|c_{r+1})$). Based on these facts and 
Proposition~\ref{prop:inc-dec}~(1), it follows that the transpositions occurring in the middle third between $\tau_{r_e}$ and $\tau_{r_l}$, which move $a_{11-r}$ from left to right, are precisely those involving elements from $X_{11-r}:=\{a_1, a_2, \ldots, a_{11-r-1}\}\cup \{b_1,\ldots ,b_r\}\cup \{c_1,\ldots ,c_r\}$. Furthermore, from Proposition~\ref{prop:inc-dec}~(1)-(2), we know that when $a_{11-r}$ enters the middle third, the elements
of $\{b_{11}, \ldots, b_{r+1}\}$ are in the initial third. 

Assertion (1) follows from the fact that $a_{11}$ goes from position 12 to position 22 by transposing with the elements of $X_{11}=\{a_1, \ldots , a_{10}\}$. Since $N^{aa}_{12}(\Pi)=2$ by Corollary~\ref{cor:xx}, the transpositions that allow $a_{11}$ to reach positions 13 and 22 must be the only
ones contributing to $N^{aa}_{12}(\Pi)$. In particular, this implies that 
$13\leq g(a_j|a_i)\leq 20$ for any $1\leq i<j\leq 10$. 

From the previous paragraph, we know that $13\leq g(a_{10}|a_j)\leq 20$ for any $1\leq j\leq 9$. Since $a_{10}$ goes from position 12 to position 22 by transposing with the elements of $X_{10}=\{a_1, \ldots , a_{9}, b_1,c_1\}$, and 
given that \(AB \prec AC \prec BC\), it follows that $a_{10}|b_1$ occurs before than $a_{10}|c_1$. Therefore, we must have $g(a_{10}|b_{1})=12$ and $g(a_{10}|c_{1})=21$. Since $N^{aa}_{13}(\Pi)=4$ by Corollary~\ref{cor:xx}, the transpositions that allow $a_{11}$ and $a_{10}$ to reach positions 14 and 21 are the only ones contributing to $N^{aa}_{13}(\Pi)$. This proves (2).  

From the previous paragraphs, we know that $14\leq g(a_{9}|a_j)\leq 19$ for any $1\leq j\leq 8$. Since $a_{9}$ goes from position 12 to position 22 by transposing with the elements of $X_{9}=\{a_1, \ldots , a_{8}, b_1,b_2,c_1,c_2\}$, and 
given that \(AB \prec AC \prec BC\), it follows that $a_{9}|b_1$ and $a_{9}|b_2$ occur before than $a_{9}|c_1$ and $a_{9}|c_2$. 
Therefore, we must have $\{g(a_{9}|b_{1}),g(a_{9}|b_{2})\}=\{12,13\}$ and 
$\{g(a_{9}|c_{1}),g(a_{9}|c_{2})\}~=~\{20,21\}$. Seeking a contradiction to assertion (3), suppose that $g(a_{9}|b_{1})~=~12$ and $g(a_{9}|b_{2})~=~13$. Then Proposition~\ref{prop:inc-dec}~(1) implies that
$g(b_1|b_2)=12$, contradicting that any transposition
contributing to $N^{bb}_{12}(\Pi)$ involves $b_{11}$. Similarly, we can
 conclude that $g(a_{9}|c_{1})=20$ and $g(a_{9}|c_{2})=21$. 
Since $N^{aa}_{14}(\Pi)=6$ by Corollary~\ref{cor:xx}, the transpositions that allow $a_{11}$, $a_{10}$ and $a_9$ to reach positions 15 and 20 are the only ones contributing to $N^{aa}_{13}(\Pi)$. This proves (3).  

Assertion (4) is an immediate consequence of (1), (2) and (3).
\end{proof}

\subsection{The vector~\ref{vec_trans2} is not realizable}
In view of Corollary~\ref{cor:xx}, to prove Theorem~\ref{thm:main}, it is enough to show that $N^{xx}_{16}(\Pi)\leq 31$. Our approach in this section can be regarded as a refinement of ideas and techniques used in~\cite{balogh2006ksets} and~\cite{avila20243symmetric}.


We recall that \(\Pi = (\pi_{0},\pi_{1},\ldots , \pi_{\binom{33}{2}})\) is defined as in Remark~\ref{rem:Pi}, and that \[\pi_{0}=(a_{11},a_{10},\ldots ,a_{1},b_{l_1},b_{l_2},\ldots ,b_{l_{11}}\-,c_{1},c_{2},\ldots ,c_{11}),\] where $b_{l_j}$ is the element of $B$ that is in the position $11+j$ of $\pi_0$. 

In light of Remark~\ref{rem:same-a-b-c}, it is sufficient to show that $N^{aa}_{16}(\Pi)\leq 31$. From now on, we will refer to the transpositions
contributing to $N^{aa}_{16}(\Pi)$ as the \dfn{halvings} of $A$. Similarly, we will say that a transposition $\tau$ occurs in the \dfn{\(k\)-center} if \(k < g(\tau)<n-k\).    

Since all monochromatic transpositions occur in the middle third, we will consider positions from the 12th to the 22nd as corresponding to the \dfn{spot} from the 1st to the 11th. For example, in permutation 
\[\pi~:=~(b_{11},\ldots ,b_{1}, a_{2}, a_{3}, a_{6}, a_{5}, a_{11}, a_{10}, a_{9}, a_{8}, a_{7}, a_{4}, a_{1}, c_{1}, \ldots , c_{11}),\] 
\(a_2\) occupies position \(\pi(a_2) = 12\), so the spot for \(a_2\), denoted \(\pospi{a_2}\), is 1, i.e. \(\pospi{a_2} = 1\); similarly \(\pospi{a_3} = 2\), and so on.
In general, if \(a_\ell\) is in position \(i\), \(12 \leq i \leq 22\), then \(\pospi{a_\ell} = \pi(a_\ell) -11 = i - 11\).
Since we will focus on the number of halvings of \(A\), we will abbreviate 
\(a_{\ell}\) as \(\ell\). For $\ell\in A$, let $[\ell]^+$ 
(respectively, $[\ell]^-$) denote the {\em number of halvings of $A$ that involve $\ell$ and move it from left to right (respectively, right to left)}.
We remark that if \(j|i\) contributes to $[j]^+$ (respectively, $[i]^-$),  then $j>i$. A general version of the following proposition was proved in~\cite{balogh2006ksets}.
\begin{prop}\label{o:[ell]+} If $j\in A$, then $[j]^+\leq \min\{2+[j]^-, j-1\}$ and $N^{aa}_{16}(\Pi)=\sum_{j=1}^{11}[j]^+=\sum_{j=1}^{11}[j]^-$.  
\end{prop}
We now introduce additional notation, in order to establish an upper bound for $[\ell]^+$ in terms of \(\pospi{1}, \ldots,\pospi{11}\). For distinct $j,l\in A=\{1,2,\ldots ,11\}$ and $\pi\in \Pi$, let us define:   
\begin{align*}
\laposg{j;\pi} &:= |\set{l > j \colon\ \pospi{l} <\pospi{j}}| &
  \text{ and  } & &
  \raposg{j;\pi} &:= |\set{l > j \colon\, \pospi{j} < \pospi{l}}|.
\end{align*}
Similarly, let
\begin{align*}
  \laposl{j;\pi} &:= |\set{l < j : \pospi{l} < \pospi{j}}| &
  \text{ and } & &
  \raposl{j;\pi} := |\set{l < j : \pospi{j} < \pospi{l}}|.
\end{align*}
If no confusion arises, we typically omit the reference to \(\pi\) in these four quantities and simply write $\pos{j}$, $\laposg{j}$, $\raposg{j}$ $\laposl{j}$ and $\raposl{j}$. 

Let \(\pi\in \Pi\) be such that $A$ occupies the middle third of $\pi$, and let
$j\in A$. A {\em {\bf go out} halving of $j$ with 
respect to $\pi$} (respectively, {\bf go in}) is a halving that contributes to $[j]^+$ and occurs after (respectively, before) $\pi$. 

We  now analyze the number of go out halvings of \(j\) based on the value of $\pospi{j}$.
If \(\pospi{j} \leq 5\), the $j$ will pass through the \(15\)-center and potentially two halvings after of $\pi$ could contribute to $[j]^+$. Additionally, if there exists  an \(l\in A\) such that \(j<l\) and \(\pospi{l} < \pospi{j}\), it is possible for \(l|j\) to occur in the $15$-center, moving $j$ from right to left and allowing for a potential increase of $[j]^+$ by 1.

If \(\pospi{j} = 6\), the transposition through which $j$ will pass  at the \(17\)-gate could contribute to $[j]^+$. Moreover, we note that each element of \(A\) greater than \(j\) that moves \(j\) to the left after of $\pi$ could also increase the value of $[j]^+$ by 1.    

If \(\pospi{j} \ge 7\), then the only way for $j$ to increase its halvings after $\pi$ is for $j$ to return to the 15-center (through 
at least $\pospi{j}-6$ transpositions of $j$ with elements of $A$ that are greater than it). Subsequently, corresponding reasoning described in previous paragraphs can be applied.     

Motivated by the previous descriptions, we define the function for the {\em expected number of go out halving of $j$} as follows.
\begin{equation}
e_{\pi(o)} (j):=
\begin{cases}
2 + \laposg{j}, & \text{if } \pos{j} \leq 5, \\
7- \pos{j} + \laposg{j}, & \text{if } \pos{j}\geq 6.
\end{cases}
\label{eq:halving-exit}
\end{equation}

We note that if \(j|i\) is a go out halving of \(j\), then $i<j$ and 
\(\pos{j}<\pos{i}\). Taking these facts into account, we define the \emph{best possible number of go out halvings of $j$ with respect to $\pi$} as:  
\begin{equation}
be_{\pi(o)}(j) \coloneqq   \min \set{ \raposl{j} , \max \set{ 0, e_{\pi(o)} (j)} }.
\label{eq:best-exit}
\end{equation}

Based in a totally analogous analysis for the number of 
{\em go in halvings of $j$}, we define the function for the {\em expected number of go in halving of $j$} as follows.

\begin{equation}
e_{\pi(i)} (j) :=
\begin{cases}
\pos{j} + \raposg{j} - 5, & \text{if } \pos{j} \leq 6, \\
2 + \raposg{j}, & \text{if } \pos{j} \geq 7.
\end{cases}
\label{eq:halving-entrance}
\end{equation}
Similarly, we define \emph{best possible number of go in halvings of $j$ with respect to $\pi$} as: 
\begin{equation}
 be_{\pi(i)}(j) := \min \set{ \laposl{j} , \max \set{ 0, e_{\pi(i)}(j)}}.
 \label{eq:best-entrance}
\end{equation}
As a consequence of the previous discussion, we have the following:
\begin{equation}
[j]^+ \leq be_{\pi(o)}(j) + be_{\pi(i)}(j)
\label{eq:best-out-in}
\end{equation}
\begin{prop}\label{p:u-single} 
If $A$ occupies the middle third of \(\pi\in \Pi\) and $j\in A$, then 
\begin{gather*}
[11]^{+} = 2, \quad [10]^{+} \leq 3, \quad [9]^{+} \leq 4, \quad [8]^{+} \leq 5, \quad [7]^{+} \leq 6,  \quad [6]^{+} \leq 5,\\
[5]^{+} \leq 4, \quad [4]^{+} \leq 3, \quad [3]^{+} \leq 2, \quad [2]^{+} \leq 1, \quad [1]^{+} = 0.
\end{gather*}
\end{prop}
\begin{proof} For $j\in \{1,\ldots ,7\}$, it follows from 
Proposition~\ref{o:[ell]+} that $[j]^+\leq j-1$. Similarly, $[11]^+=2$ was established in Proposition~\ref{eta_12+r}~(1). The remaining inequalities can be easily derived from Inequality~\ref{eq:best-out-in}.   
\end{proof}
 
For the remainder of the paper, we will use $[j|i]$ to mean that $j|i$ is a halving contributing to $[j]^+$. Similarly, for brevity, we define $h(\Pi):=N^{aa}_{16}(\Pi)=\sum_{j=1}^{11}[j]^+$ and assume that \(\pi\) is a permutation of $\Pi$ in which $A$ occupies the middle third.

\begin{cor}\label{cor:11-8<=13}
If $A$ occupies the middle third of \(\pi\in \Pi\), then 
\[[11]^+ + [10]^+ + [9]^- + [8]^+ \leq 13.\]
\end{cor}
\begin{proof}
From Proposition~\ref{p:u-single}, it follows that \([11]^+ + [10]^+ + [9]^- + [8]^+ \leq 14\), and equality only holds if each \(j \in \set{8,9,10,11}\) reaches its maximum. If \([10]^+ = 2\), there is nothing to prove. 
So, we may assume that \([10]^+ = 3\). 
From Proposition~\ref{o:[ell]+}, it follows that \([10]^- = 1\), which implies that \([11|10]\) holds.
As \([11]^+ = 2\) and \([11|10]\) holds, then \(11\) is involved in a halving with at most one of 8 or 9. 
Thus, \([9]^+ + [8]^+ \leq 8\), so \([11]^+ + [10]^+ + [9]^- + [8]^+ \leq 13\).
\end{proof}

\begin{cor}\label{c:>=14} 
If $[7]^+ + [8]^+ + [9]^+ + [10]^+\leq 14$, then
$h(\Pi)\leq 31$. \hfill \qedsymbol
\end{cor}

\begin{prop}\label{p:at-most-one} If \(h(\Pi)\geq 32\), then 
$12\leq [7]^++[8]^++[9]^+\leq 13$. If
additionally \([10]^{+} = 3\), then at most one of $[7]^+$, $[8]^+$,$[9]^+$ reaches the upper bound given in Proposition~\ref{p:u-single}. 
\end{prop} 
\begin{proof} Given \(h(\Pi)\geq 32\), $[10]^+\leq 3$ and Corollary~\ref{c:>=14}, it follows that 
$[7]^++[8]^++[9]^+\geq 12$. 

Let $i,j\in \{7,8,9\}$ with $j>i$. Suppose that 
$[i]^+=(11-i)+2$ and $[j]^+=(11-j)+2$ (reaching their maximum possible values). This implies that
$[11|l]$ and $[10|l]$ must hold for each $l\in \{i,j\}$.
Consequently, $[11|10]$ cannot hold, which means that $[10|i]$ and $[10|j]$ are the only halvings contributing to $[10]^+$. Then
if $k\in \{7,8,9\}\setminus \{i,j\}$, we must have
$[k]^-\leq (11-k)-2$, and therefore
$[k]^+\leq 2+((11-k)-2)=11-k$.
Then $[7]^++[8]^++[9]^+\leq (11-7)+(11-8)+(11-9)+4=13$.

For the second assertion, we only show the case $[9]^+=4$, because the other cases are entirely analogous.  Since $[10]^+=3$, it follows that $[11|10]$ must hold. Similarly, $[9]^+=4$ implies that $[11|9]$ is true. Given that $[11]^+=2$, neither $[11|8]$ nor $[11|7]$ can hold. These facts imply that $[8]^-\leq 2$ and $[7]^-\leq 3$. Assertion (A) follows by applying $[j]^+\leq 2+[j]^-$ to $j=7,8$.
\end{proof}

\begin{cor}\label{c:general-bounds} If \(h(\Pi)\geq 32\), then $4\leq [10]^++[11]^+\leq 5$ and 
$14\leq [6]^++[5]^++\ldots +[1]^+\leq 15$. 
\end{cor} 
\begin{proof}
The upper bounds follow directly from 
Proposition~\ref{p:u-single}. The lower bounds follow from
the upper bounds and the upper bound given in Proposition~\ref{p:at-most-one}. 
\end{proof}

\begin{prop}\label{p:7,8;4,3} If \(h(\Pi)\geq 32\), then \([8]^{+} \geq 3\) and $[7]^+\geq 4$.
\end{prop} 
\begin{proof} To seek a contradiction, suppose that $[7]^+\leq 3$. If $[10]^+=2$, then 
Proposition~\ref{p:u-single} implies $h(\Pi)\leq 2+2+4+5+3+15=31$, a contradiction. Therefore, we can assume that $[10]^+=3$. From this and Proposition~\ref{p:at-most-one}, it follows that $[8]^++[9]^+\leq 8$, which implies $h(\Pi)\leq 2+3+8+3+15=31$, a contradiction. This proves $[7]^+\geq 3$.

Similarly, suppose that $[8]^+\leq 2$. If $[10]^+=2$, then Proposition~\ref{p:u-single} implies $h(\Pi)\leq 2+2+4+2+6+15=31$, a contradiction. We then assume that $[10]^+=3$. From this and Proposition~\ref{p:at-most-one}, it follows that $[7]^++[9]^+\leq 9$, leading to $h(\Pi)\leq 2+3+9+2+15=31$. This contradiction proves that
$[8]^+\geq 3$.    
\end{proof}

In any allowable sequence \(\Pi\), there is a permutation \(\pi\) such that \(\pospi{11}=5\), since, by assumption, 11 (\(=a_{11}\)) starts in the initial third and ends in the final third. In the following propositions, we will assume this fact without loss of generality.

\begin{prop}\label{p:15-center} If \(\pospi{11}=5\), $r\in \{7, 8,9,10\}$ and \([r]^{+} =11-r+2\), then \(r\) is in the 15-center of $\pi$.
\end{prop} 
\begin{proof}
Let $r\in \{7,8,9,10\}$. From \([r]^{+} =11-r+2\), it follows that $l|r$ must be a halving for any  $l\in \{r+1,\ldots,11\}$. Since \(\pospi{r}\leq 4\) implies that $11|r$ cannot be a halving, it follows that \(\pospi{r}\geq 5\). Similarly, since  \(\pospi{r}\geq 8\) implies that $m|r$ cannot be a halving for some $m\in \{r+1,\ldots,11\}$, the  assertion follows. 
\end{proof}
\begin{prop}\label{p:[r]+=r-1]} If $\pospi{11}=5$, $r\in \{4,5,6,7\}$ and \([r]^{+} =r-1\), then \(r\) is in the  ($r+8$)-center of $\pi$.
\end{prop}
\begin{proof} Let $r\in \{4,5,6,7\}$ and $l:=\pospi{r}$. 
To seek a contradiction, suppose that $l\leq r-3$ or $l\geq 15-r$. 
We suppose first that $l\leq r-3$. We note that
$\laposl{\pi;r} +\laposg{\pi;r}=l-1$. It is easy to see that 
$be_{\pi(o)}(r)\leq 2+\laposg{\pi;r} $ and 
$be_{\pi(i)}(r)\leq \laposl{\pi;r}$. From these and~\ref{eq:best-out-in} we obtain $[r]^+\leq  (2+\laposg{\pi;r}) + \laposl{\pi;r}\leq 
(l-1)+2=l+1\leq r-2$, a contradiction. 

Suppose now that $l\geq 15-r$. We note that 
$\raposl{\pi;r} +\raposg{\pi;r}=11-l$. It is easy to see that 
$be_{\pi(o)}(r)\leq \raposl{\pi;r}$ and 
$be_{\pi(i)}(r)\leq 2+\raposg{\pi;r}$. From these and~\ref{eq:best-out-in} we obtain $[r]^+\leq \raposl{\pi;r}+(2+\raposg{\pi;r})\leq 13-l\leq 13-(15-r)=r-2$, a contradiction. 
\end{proof}

\begin{prop}\label{p:8,9,10}
If $h(\Pi)\geq 32$, $\pospi{11}=5$ and $j\in\{8,9,10\}$, then $4\leq \pospi{j}\leq 8$. 
\end{prop}

\begin{proof} Let $p\in\{8,9,10\}$ be such that 
$\pospi{p}:=\min\{\pospi{8},\pospi{9},\pospi{10}\}$. 
To seek a contradiction, suppose first that $\pospi{p}\leq 3$. Let $T$ be the set of  transpositions that move $a_p$ from (left to right) position $\pospi{p}$ to position 4, and then from position $8$ to position $11$. We have 
$|T|\geq (4-\pospi{p})+3\geq 4$. Since all transpositions in $T$ occur outside of the $14$-center, Proposition~\ref{eta_12+r}~(4) implies that none of these can be monochromatic. Thus, each transposition of $T$ must be of type $ac$. On the other hand, since 
$g(a_p|c_{11-p+1})=22$ by Proposition~\ref{prop:inc-dec}~(3), we can deduce from Proposition~\ref{prop:inc-dec}~(1) that $|T|\leq 11-p\leq 3$, contradicting that $|T|\geq 4$.   

Similarly, let $q\in\{8,9,10\}$ be such that 
$\pospi{q}:=\max\{\pospi{8},\pospi{9},\pospi{10}\}$. 
To seek a contradiction, suppose now that 
$\pospi{q}\geq 9$. Let $T$ be the set of  transpositions that move $a_q$ from (left to right) position $1$ to position 4, and then from position $8$ to position $\pospi{q}$.  We note that each
transposition of $T$ occur outside of the $14$-center.

Since $g(a_{q}|b_{12-q})=11$ by 
Proposition~\ref{prop:inc-dec}~(2) and $|T|\geq 3+(\pospi{q}-8)\geq 4$, then at most 3 (respectively, at least 1) transpositions of $T$ are of type $ab$ (respectively, of type $aa$). 
This fact and Proposition~\ref{eta_12+r}~(4), imply that $q\in\{9,10\}$. However, this last implies that some of $11|9,~11|10$ or $10|9$ contributes to 
$N^{aa}_{12}(\Pi)+N^{aa}_{13}(\Pi)+N^{aa}_{14}(\Pi)$, which contradicts Proposition~\ref{eta_12+r}~(4). 
\end{proof}

\begin{prop}\label{p:[6]=5,[7]<5} If 
$h(\Pi)\geq 32$ and $[6]^{+}=5$, then $[7]^{+} = 4$.
\end{prop}

\begin{proof} 
Suppose that $[6]^{+}=5$.
Let $\pi \in \Pi$ be such that $\pospi{11}=5$.
By Proposition~\ref{p:8,9,10}, we know that $8,~9$ and $10$ are in the $14$-center of $\pi$. 
Similarly, as $[6]^+=5$, Proposition~\ref{p:[r]+=r-1]} implies that $6$ must be in the $14-$center of $\pi$.  
Thus, the $14-$center of $\pi$ is occupied by  $6, 8, 9, 10, 11$. 
Clearly, if $l:=\pospi{7}$, then either $l\leq 3$ or $l\geq 9$.

If $l\leq 3$, then $\laposl{\pi;7} +\laposg{\pi;7}=l-1$. It is easy to see that 
$be_{\pi(o)}(7)\leq 2+\laposg{\pi;r} $ and 
$be_{\pi(i)}(7)\leq \laposl{\pi;7}$. 
From these and Equation~\ref{eq:best-out-in} we obtain $[7]^+\leq  (2+\laposg{\pi;7}) + \laposl{\pi;7}\leq 
(l-1)+2=l+1\leq 4$. 

If $l\geq 9$, then $\raposl{\pi;7} +\raposg{\pi;7}=11-l$. It is easy to see that 
$be_{\pi(o)}(7)\leq \raposl{\pi;7}$ and 
$be_{\pi(i)}(7)\leq 2+\raposg{\pi;7}$. 
From these and Equation~\ref{eq:best-out-in} we obtain $[7]^+\leq \raposl{\pi;7}+(2+\raposg{\pi;7})\leq 13-l\leq 4$.  

To seek a contradiction, suppose that \([7]^+ < 4\). 
We know that either \([10]^+ = 3\) or \([10]^+ = 2\).
If \([10]^+ =2\), and each \([j]^+\), for \(j \in \set{1, 2, \ldots, 11} \setminus \set{7,10}\), reaches its maximum given by Proposition~\ref{p:u-single}, a simple calculation yields that \(h(\Pi) \leq 31\), a contradiction.
If \([10]^+ = 3\), by Proposition~\ref{p:at-most-one}, at most one of \([7]^+, [8]^+\), and \([9]^+\) reaches its maximum given by Proposition~\ref{p:u-single}, then \([8]^+ + [9]^+ \leq 8\). A simple calculation yields that \(h(\Pi) \leq 31\), a contradiction.
Therefore, \([7]^+ =4\).
\end{proof} 

Let $r\in \{3,4,5,6\}$. If $[i]^+=i-1$ for each 
$i\in \{1,\ldots ,r\}$, we will say that $K_r$ \dfn{holds}.

\begin{prop}\label{p:pre-last} Let $\pi\in \Pi$ be such that 
\(\pospi{11} = 5\) and \(h(\Pi) \geq 32\). Then the following statements hold:
\begin{itemize}
    \item[(1)] If 6 lies in the 14-center of $\pi$, then \(\{\pospi{5}, \pospi{7}\}=\set{3,9}\) and \(\pospi{4} \in \set{2,10}\).
    \item[(2)] If 6 does not lie in the 14-center of $\pi$, then \(\{\pospi{5}, \pospi{6}\}=\set{3,9} \) and  \(\pospi{4} \in \set{2,10}\).
\end{itemize} 
\end{prop}
\begin{proof} 
From Proposition~\ref{p:8,9,10}, we know that 8, 9, 10 and 11 are in some of the 5 spots of the 14-center of $\pi$. So, only one of 1, 2, \ldots, 7 could be in the 14-center of $\pi$. 

Suppose first that 6 lies in the 14-center of $\pi$. Since \(h(\Pi) \geq 32\), Proposition~\ref{p:7,8;4,3} yields 
\([7]^+ \geq 4\). Since either \(\pospi{7} \leq 3\) or \(\pospi{7} \geq 9\), then \([7]^+ \leq 4\). Thus \([7]^+ = 4\), and so we must have that \(\pospi{7} \in \set{3, 9}\).
By Corollary~\ref{cor:11-8<=13}, \([11]^+ + [10]^+ + [9]^+ + [8]^+ \leq 13\). Consequently,  \([6]^+ + \cdots + [1]^+ = 15\), i. e. \(K_6\) holds. Then \([r]^+ = r-1\) for each \(r \in \set{6, 5,4,3,2,1}\).
This last and Proposition~\ref{p:[r]+=r-1]} imply that
\(\{\pospi{5}, \pospi{7}\}=\set{3,9}\) and \(\pospi{4} \in \set{2,10}\).
This proves~(1). 

Suppose now that 6 does not lie in the 14-center of $\pi$.
Then \(\pospi{6} \leq 3\) or \(\pospi{6} \geq 9\). 
By Inequality~\ref{eq:best-out-in}, we have \([6]^+ \leq 4\). Moreover, since \(h(\Pi) \geq 32\), Corollary~\ref{c:general-bounds} implies \(14 \leq [6]^+ + \cdots + [1]^+ \leq 15\). Hence \([6]^+ =4\), which forces \(\pospi{6} \in \set{3,9}\). 
In addition, \([6]^+ + \cdots + [1]^+ =14\) and \([r]^+ = r-1\), for \(r \in \set{5,4,3,2,1}\).
If \(\pospi{7} \leq 2\) or \(\pospi{7} \geq 10\), then \([7]^+ \leq 3\), which contradicts Proposition~\ref{p:7,8;4,3}.
If \(\pospi{7} \in \set{9,3}\), then \([7]^+ \leq 4\), in fact we must have \([7]^+ = 4\). 
As \([11]^+ + [10]^+ + [9]^+ + [8]^+ \leq 13\) and \([6]^+ + \cdots + [1]^+ = 14\), then \([11]^+ + \cdots + [1]^+ \leq 31\), which is a contradiction. Thus, 7 must be in the 14-center of $\pi$.
Because \([r]^+ = r-1\) for each \(r \in \set{5,4,3,2,1}\), then Proposition~\ref{p:[r]+=r-1]} implies
\(\{\pospi{5}, \pospi{6}\}=\set{3,9} \) and  \(\pospi{4} \in \set{2,10}\).
\end{proof}

The following observation is implicit in the proof of Proposition~\ref{p:pre-last}.  

\begin{obs}\label{obs:K5-holds}
If $h(\Pi)\geq 32$, then \(K_5\) holds; that is, \([r]^+ = r-1\) 
for all \(r \in \set{5,4,3,2,1}\).
\end{obs}

\begin{lem}
If $h(\Pi)\geq 32$, then $\Pi$ is not realizable.
\end{lem}
\begin{proof} Let $\pi\in \Pi$ be such that $\pospi{11}=5$. 
From Proposition~\ref{p:pre-last} we know that the elements in \(\{11,10,9,8,7,6,5\}\) occupy the 13-center of $\pi$, \(\pospi{5} \in \{3,9\}\)  and that \(\pospi{4} \in \{2,10\}\). Therefore, four cases arise:

\begin{itemize}[leftmargin=2cm]
    \item[Case (A).] \(\pospi{4}=2\) and \(\pospi{5}=3\) \label{case:sp(4)-2,sp(5)-3}
    \item[Case (B).] \(\pospi{4}=2\) and \(\pospi{5}=9\) \label{case:sp(4)-2,sp(5)-9}
    \item[Case (C).] \(\pospi{4}=10\) and \(\pospi{5}=9\)  \label{case:sp(4)-10,sp(5)-9}
    \item[Case (D).] \(\pospi{4}=10\) and \(\pospi{5}=3\) \label{case:sp(4)-10,sp(5)-3}
\end{itemize}

We shall prove only the Cases~(A) and~(B), because the Cases~(C) and~(D) are analogous, respectively, and they occur when \(B\) is the middle third and \(B\) is interchanging with \(C\) (which is in the initial third). 

We recall that, when analyzing permutations with \(A\) in the middle third, the element \(a_\ell\) is denoted simply by \(\ell\). So, when we say \(\pospi{5} = 3\), it actually means \(\pospi{a_5} = 3\). In order to avoid confusions, we will take back the usual notation.

\begin{claim}
Case~(A) is not realizable.
\end{claim}

Let
\[\pi=(B,x_{1},a_{4},a_{5},a_{l},a_{11},a_{m},a_{n},a_{p},a_{q},x_{2},x_{3},C),\]
where \(\set{x_{1},x_{2},x_{3}} = \set{a_{1},a_{2},a_{3}}\) and \(\set{l,m,n,p,q} =  \set{10,9,8,7,6}\). 
After \(\{a_{6},\ldots ,a_{11}\}\) leaves the middle third in descending order, as we know, we have
\[(B,y_{3},y_{2},x_{1},a_{4},a_{5},x_{2},x_{3},y_{1},c_{4},c_{5},c_{6},X),\]
with \(\set{y_{1},y_{2},y_{3}} = \set{c_{1},c_{2},c_{3}}\).
Afterwards \(a_{5}\) leaves the middle third we have:
\[(B,y_{3},y_{2},x_{1}, a_4, \colorbox[RGB]{225,225,225}{\(x'_{2}, x'_{3}, y_{1},\)} c_4, c_{5},c_{6},c_{7},X),\]
where \(\set{x'_{2},x'_{3}} = \set{x_{2},x_{3}}\). 

Then,  \(a_{4} | x_{2}^{'}\) or \(c_{4} | y_{1}\) are not halvings, which contradicts that \(K_5\) holds (by symmetry) for each of \(A, B, \) and \(C\).

\begin{claim}
Case~(B) is not realizable.
\end{claim}

Let
\[\pi=(B,x_{1},a_{4},a_{k}, a_{l},a_{11},a_{m},a_{n},a_{p},a_{5},x_{2},x_{3},C),\]
where \(\set{x_{1},x_{2},x_{3}} = \set{a_{1},a_{2},a_{3}}\) and \(\set{k,l,m,n,p} = \set{10,9,8,7,6}\).
After \(\{a_{6},\ldots ,a_{11}\}\) leaves the middle third in descending order we have
\[(B,y_3,y_{2},x_{1},a_{4},a_{5},x_{2},x_{3},y_{1},c_{4},c_{5},c_{6},X),\]
with \(\set{y_{1},y_{2},y_{3}} = \set{c_{1},c_{2},c_{3}}\).
Afterwards \(a_{5}\) leaves the middle third we have
\[(B,y_{3},y_{2},x_{1},a_{4},\colorbox[RGB]{225,225,225}{\(x'_{2},x'_{3},y_{1},\)} c_{4},c_{5},c_{6},c_{7},X),\]
where \(\set{x'_{2},x'_{3}} = \set{x_{2},x_{3}}\). 
Then,  \(a_{4} | x_{2}^{'}\) or \(c_{4} | y_{1}\) are not halvings, which contradicts that \(K_5\) holds (by symmetry) for each of \(A, B, \) and \(C\).
\end{proof}

\section{Acknowledgments}

V\'ictor H. G\'omez Mart\'inez is a doctoral student of the Programa de Doctorado en Ciencias Aplicadas at Universidad Autónoma de San Luis Potos\'i (UASLP) and has received a fellowship (Grant No. 370207) from Secretar\'ia de Ciencia, Humanidades, Tecnolog\'ia e Innovaci\'on \makebox{({SECIHTI})}, formerly  Consejo Nacional de Humanidades, Ciencias y Tecnolog\'ias \makebox{(CONAHCYT)}.


\bibliographystyle{plain}
\bibliography{3-sym-pseudo-K33-biblio.bib}

\end{document}